\documentclass[english]{article}
\usepackage[T1]{fontenc}
\usepackage[latin9]{inputenc}
\usepackage{geometry}
\usepackage{amsmath}
\usepackage{amsthm}
\usepackage{amssymb}
\usepackage{stmaryrd}
\usepackage{rotfloat}
\usepackage{setspace}
\makeatletter

\providecommand{\tabularnewline}{\\}

\theoremstyle{plain}
\newtheorem{thm}{\protect\theoremname}[section]
\theoremstyle{plain}
\newtheorem{cor}[thm]{\protect\corollaryname}
\theoremstyle{plain}
\newtheorem{lem}[thm]{\protect\lemmaname}

\usepackage[all]{xy}
\xyoption{rotate}
\usepackage{tocbibind}
\usepackage{indentfirst}
\usepackage{setspace}
\usepackage{bm}

\makeatother

\usepackage{babel}
\providecommand{\corollaryname}{Corollary}
\providecommand{\lemmaname}{Lemma}
\providecommand{\theoremname}{Theorem}

\begin{document}
\title{Integral Motivic Cohomology of $BSO_{4}$}
\author{Alexander Port}
\maketitle
\begin{abstract}
Motivic cohomology is powerful tool in algebraic geometry with associated
realization maps giving important information about the relations
between cohomological invariants of schemes and their classifying
spaces. The problem of computing general cohomological invariants
of these classifying spaces is ongoing. Most relevant to this paper
is (1) Totaro's construction of the Chow ring of a classifying space
in general and his use of this to study symmetric groups in \cite{Tot99},
(2) Guillot's similar examination for the Lie groups $G_{2}$ and
$Spin(7)$ in \cite{Gui07}, (3) Field's computation of the Chow ring
of $BSO(2n,\mathbb{C})$ in \cite{Fie12}, and (4) Yagita's work on
the $\mathbb{Z}_{2}$-motivic cohomology of $BSO_{4}$ and $BG_{2}$
in \cite{Yag10}.

The work presented in this paper covers the computation of the motivic
cohomology of $BSO_{4}$ with integral coefficients. The primary approach
draws on methods laid out by Guillot and Yagita (\cite{Gui07}, \cite{Yag10}).
These results lay the groundwork for future work, most immediately
the analogous computation for $BG_{2}$ (\cite{Por21}).
\end{abstract}

\section{Notation and Conventions}
\begin{itemize}
\item Maps:
\begin{itemize}
\item The short exact sequence $0\rightarrow\mathbb{Z}\xrightarrow{\times2}\mathbb{Z}\rightarrow\mathbb{Z}_{2}\rightarrow0$
induces the long exact sequence
\end{itemize}
\[
...\rightarrow H^{*}(BSO_{4};\mathbb{Z})\xrightarrow{\times2}H^{*}(BSO_{4};\mathbb{Z})\xrightarrow{\mu}H^{*}(BSO_{4};\mathbb{Z}_{2})\xrightarrow{\tilde{\beta}}H^{*+1}(BSO_{4};\mathbb{Z})\rightarrow...
\]

\begin{itemize}
\item The short exact sequence $0\rightarrow\mathbb{Z}_{2}\xrightarrow{\times2}\mathbb{Z}_{4}\rightarrow\mathbb{Z}_{2}\rightarrow0$
induces a similar long exact sequence
\end{itemize}
\[
...\rightarrow H^{*}(BSO_{4};\mathbb{Z}_{2})\xrightarrow{\times2}H^{*}(BSO_{4};\mathbb{Z}_{4})\rightarrow H^{*}(BSO_{4};\mathbb{Z}_{2})\xrightarrow{\beta}H^{*+1}(BSO_{4};\mathbb{Z}_{2})\rightarrow...
\]

\begin{itemize}
\item Note that $\mu$ is the map induced by reducing coefficients mod $2$,
$\tilde{\beta}$ is the \emph{integral Bockstein homomorphism}, and
$\beta$ is the \emph{Bockstein homomorphism} (where $\beta=\mu\circ\tilde{\beta}$).
\item The Milnor operations $Q_{i}:H^{*,*'}(X;\mathbb{Z}_{p})\rightarrow H^{*+2p^{i}-1,*'+p^{i}-1}(X;\mathbb{Z}_{p})$
are defined inductively as $Q_{0}=\beta$ and $Q_{i+1}=[Q_{i},P^{p^{i}}]$
where $P^{k}:H^{*,*'}(X;\mathbb{Z}_{p})\rightarrow H^{*+2k(p-1),*'+k(p-1)}(X;\mathbb{Z}_{p})$
is a reduced $p^{th}$-power operation.
\end{itemize}
\item Generators:
\begin{itemize}
\item In the context of ordinary cohomology and Chow rings, the subscript
of a generator is usually its degree. For motivic cohomology, a subscript
of $n$ means the degree is $(2n,n)$. The one exception to this is
$H^{*,*'}((\mathbb{G}_{m})^{n};\mathbb{Z})=H^{*,*'}(Spec(\mathbb{C});\mathbb{Z})[t_{1},...,t_{n}]$
where $deg(t_{i})=(2,1)$ for all $i$ (\cite{Fie12}, \cite{Gui07}).
\item Let $d_{m}\in H^{2m}(BSO_{n};\mathbb{Z})$ be the $m^{th}$ Chern
class of the standard representation of $SO_{n}$ (\cite{Gui07}).
\item Let $\zeta$ be the universal $n$-plane bundle over $BSO_{n}$, let
$w_{m}\in H^{m}(BSO_{n};\mathbb{Z}_{2})$ be the $m^{th}$ Stiefel-Whitney
class of $\zeta$, and let $p_{m}\in H^{4m}(BSO_{n};\mathbb{Z})$
be the $m^{th}$ Pontryagin class of the complexification $\zeta\otimes\mathbb{C}$
(\cite{Bro82}).
\item $\tau$ refers to a generator of $H^{*,*'}(Spec(\mathbb{C});\mathbb{Z}_{p})$
with $deg(\tau)=(0,1)$ and corresponds to a choice of generator of
$H^{0,1}(Spec(\mathbb{C});\mathbb{Z}_{p})\cong\mathbb{Z}_{p}$.
\end{itemize}
\item Shorthand:
\begin{itemize}
\item We use the shorthand $Q(n)$ to refer to the exterior algebra $\Lambda_{\mathbb{Z}_{2}}(Q_{0},...,Q_{n})$
(\cite{Yag10}).
\item If $I=\{i_{1},...,i_{l}\}$ we denote the products $w(2I)=w_{2i_{1}}...w_{2i_{l}}$
and $p(I)=p_{i_{1}}...p_{i_{l}}$.
\item If $S$ is a set and $R$ is a commutative ring then we use $R\{S\}$
to denote the free $R$-module with $S$ as its basis.
\item If $A$ is an abelian group then $_{l}A=\{a\in A:la=0_{A}\}$ is the
$l$-torsion subgroup of $A$.
\end{itemize}
\end{itemize}

\section{Known Computations}

\subsection{Classical $BSO_{4}$ Computations}

The computational strategy used for computing $H^{*,*'}(BSO_{4};\mathbb{Z})$
is to take lifts of classical generators (against the realization
map $t:H^{*,*'}(BSO_{4};A)\rightarrow H^{*}(BSO_{4};A)$ for $A=\mathbb{Z}$
or $\mathbb{Z}_{p}$) and to search for additional generators using
the universal coefficient theorem. As such, the following represents
a summary of the needed classical results regarding the cohomology
and Chow ring of $BSO_{4}$:
\begin{thm}
(\cite{Bro82} 1.3) The mod $2$ cohomology ring of $BSO_{n}$ is
given by

\[
H^{*}(BSO_{n};\mathbb{Z}_{2})=\mathbb{Z}_{2}[w_{2},w_{3},...,w_{2n}]
\]
\end{thm}

\medskip{}

\begin{thm}
(\cite{Bro82} 1.5) The integral cohomology ring of $BSO_{n}$ is
$H^{*}(BSO_{n};\mathbb{Z})=\bar{R}_{n}/\bar{I}_{n}$ where

\[
\bar{R}_{n}=\mathbb{Z}[p_{1},...,p_{\left\lfloor \frac{n-1}{2}\right\rfloor },X_{n},\tilde{\beta}w(2I):1\leq i_{1}<...<i_{l}\leq\left\lfloor \frac{n-1}{2}\right\rfloor ]
\]
 and $\bar{I}_{n}$ is the ideal generated by the following relations:
(1) $2\tilde{\beta}w(I)=0$, (2) $\tilde{\beta}w(2I)\tilde{\beta}w(2J)$
is equal to

\[
\sum_{k\in I}\tilde{\beta}w_{2k}p((I-\{k\})\cap J)\tilde{\beta}w(2((I-\{k\})\cup J-(I-\{k\})\cap J))
\]
 (3) if $n=2k+1$ is odd then $X_{n}=\tilde{\beta}w_{2k}$, and (4)
if $n=2k$ is even then $X_{n}^{2}=p_{k}$. Furthermore $\mu(p_{i})=w_{2i}^{2}$
for all $i$.
\end{thm}

\medskip{}

\begin{thm}
(\cite{Fie12} 1) The Chow ring of $BSO_{n}$ is given by

\[
CH^{*}(BSO_{n})=\mathbb{Z}[d_{2},...,d_{2n},y_{n}]/(2d_{odd},y_{n}d_{odd},y_{n}^{2}+(-1)^{n}2^{2n-2}d_{2n})
\]
 where $y_{n}$ is the class which maps to $2^{n-1}$ times the Euler
class of $\zeta_{n}$ under the cycle class map.
\end{thm}

\noindent Note that when $n=2k$ is even, the above class $X_{n}$
is the Euler class $e=e(\zeta\otimes\mathbb{C})$ as is satisfies
the needed property of $e\smile e=p_{k}$. As such, we will simply
replace $X_{n}$ with $\sqrt{p_{k}}$ moving forward.
\begin{cor}
(\cite{Por21} 5.3.4)\label{cor: classical rings of BSO4}Taking $n=4$
with the above theorems yields
\end{cor}

\noindent 
\begin{equation}
\begin{array}{c}
H^{*}(BSO_{4};\mathbb{Z}_{2})=\mathbb{Z}_{2}[w_{2},w_{3},w_{4}]\\
H^{*}(BSO_{4};\mathbb{Z})=\mathbb{Z}[\tilde{\beta}w_{2},p_{1},\sqrt{p_{2}}]/(2\tilde{\beta}w_{2})\\
CH^{*}(BSO_{4})=\mathbb{Z}[d_{2},d_{3},d_{4},y_{2}]/(2d_{3},y_{2}d_{3},y_{2}^{2}-4d_{4})
\end{array}
\end{equation}

\subsection{Re-framing Yagita's $H^{*,*'}(BSO_{4};\mathbb{Z}_{2})$}

With this in place we can now move into the motivic setting. Recalling
that $H^{2n,n}(X;A)=CH^{n}(X)\otimes A$, we can take $d_{2},d_{3},d_{4}$
to be motivic cohomology classes where

\begin{equation}
\begin{array}{ccc}
deg(d_{2})=(4,2) & deg(d_{3})=(6,3) & deg(d_{4})=(8,4)\end{array}
\end{equation}

\noindent However, to avoid confusion we denote the motivic Stiefel-Whitney
classes as $w_{i}^{alg}$ (i.e. the ``algebraic'' lifts of their
classical counterparts). As seen by Corollary 2.2 of \cite{HN18}:

\begin{equation}
\begin{array}{cc}
deg(w_{2i}^{alg})=(2i,i+1) & deg(w_{2i+1}^{alg})=(2i+1,i+1)\end{array}
\end{equation}

\noindent With these generators and conventions in place we can consider
Yagita's original statement of the mod $2$ motivic cohomology of
$BSO_{4}$:
\begin{thm}
(\cite{Yag10} 9.4) The motivic cohomology $H^{*,*'}(BSO_{4};\mathbb{Z}_{2})$
is given by

\[
\begin{array}{c}
H^{*,*'}(BSO_{4};\mathbb{Z}_{2})=\mathbb{Z}_{2}[\mu(d_{2}),\mu(d_{4})]\{y_{0,2}\}\oplus\mathbb{Z}_{2}[\tau,\mu(d_{2})]\otimes(\mathbb{Z}_{2}[\mu(d_{4})]\{1\}\oplus\mathbb{Z}_{2}[\mu(d_{3})]\otimes Q(1)\{w_{2}^{alg}\}\\
\oplus\mathbb{Z}_{2}[\mu(d_{4})]\otimes(\mathbb{Z}_{2}[\mu(d_{3})]Q(2)-\mathbb{Z}_{2}\{1\})\{a\})
\end{array}
\]

\noindent where (1) $a$ is a virtual element corresponding to element
that in $H^{0}(BSO_{4};H^{3}(B\mathbb{Z}_{2};\mathbb{Z}_{2}))$ and
is taken be have degree $(3,3)$, and (2) the degrees of the generators
are given by $deg(\tau)=(0,1)$, $deg(y_{0,2})=(4,2)$, and $deg(\mu(d_{i}))=(2i,i)$.
Furthermore, realizations involving $a$ are given as follows:

\[
\begin{array}{cccc}
t_{2}(\mu(d_{2})a)=w_{2}w_{3}w_{4} & t_{2}(Q_{0}a)=w_{4} & t_{2}(Q_{1}a)=w_{2}w_{4} & t_{2}(Q_{2}a)=w_{2}w_{4}^{2}\end{array}
\]
\end{thm}

\noindent The following theorem will re-frame Yagita's above generators
in a different context which will be useful in the following section.
\begin{thm}
(\cite{Por21} 7.3.3) The motivic cohomology ring $H^{*,*'}(BSO_{4};\mathbb{Z}_{2})$
is given by
\[
\mathbb{Z}_{2}[\tau,y_{0,2},w_{2}^{alg},w_{3}^{alg},w_{4}^{alg},\tau^{-2}(w_{2}^{alg})^{2},\tau^{-1}(w_{3}^{alg})^{2},\tau^{-2}(w_{4}^{alg})^{2},\tau^{-1}w_{2}^{alg}w_{3}^{alg},\tau^{-1}w_{2}^{alg}w_{4}^{alg},\tau^{-1}w_{3}^{alg}w_{4}^{alg}]/I
\]

\noindent where degrees are given by

\noindent 
\[
\begin{array}{ccccc}
deg(\tau)=(0,1) & deg(y_{0,2})=(4,2) & deg(w_{2}^{alg})=(2,2) & deg(w_{3}^{alg})=(3,2) & deg(w_{4}^{alg})=(4,3)\end{array}
\]

\noindent and $I$ is the ideal generated by the following relations:

\[
\begin{array}{c}
\begin{array}{ccccc}
\tau y_{0,2}=0 & y_{0,2}^{2}=0 & y_{0,2}w_{2}^{alg}=0 & y_{0,2}w_{3}^{alg}=0 & y_{0,2}w_{4}^{alg}=0\end{array}\\
\begin{array}{cccc}
y_{0,2}\tau^{-1}(w_{3}^{alg})^{2}=0 & y_{0,2}\tau^{-1}w_{2}^{alg}w_{3}^{alg}=0 & y_{0,2}\tau^{-1}w_{2}^{alg}w_{4}^{alg}=0 & y_{0,2}\tau^{-1}w_{3}^{alg}w_{4}^{alg}=0\end{array}
\end{array}
\]
\end{thm}

\noindent Finally, Harada and Nakada expand on Yagita's result and
show that $Ker(t_{2})=\mathbb{Z}_{2}[\mu(d_{2}),\mu(d_{4})]\{y_{0,2}\}$
(\cite{HN18} 0.1); translating this into the new notation we get
that

\begin{equation}
Ker(t_{2})=\mathbb{Z}_{2}[\tau^{-2}(w_{2}^{alg})^{2},\tau^{-2}(w_{4}^{alg})^{2}]\{y_{0,2}\}
\end{equation}

\section{Integral Motivic Cohomology of $BSO_{4}$ }

The goal here is to compute the integral motivic cohomology of $BSO_{4}$
using the universal coefficient theorem (UCT) as presented on pg.
27 of ``Lecture Notes on Motivic Cohomology'' by Carlo Mazza, Vladimir
Voevodsky and Charles Weibel as well its classical analog (\cite{MVW06},
see \cite{nLab} universal coefficient theorem for classical version
used here). Here we obtain a short exact sequence on cohomology:

\[
0\rightarrow H^{*}(BSO_{4};\mathbb{Z})\otimes\mathbb{Z}_{2}\xrightarrow{\mu_{C}}H^{*}(BSO_{4};\mathbb{Z}_{2})\xrightarrow{\tilde{\beta}_{C}}{}_{2}H^{*+1}(BSO_{4};\mathbb{Z})\rightarrow0
\]

\noindent These maps $\mu_{C}$ and $\tilde{\beta}_{C}$ are closely
related to $\mu$ and $\tilde{\beta}$ in that $\mu(x)=\mu_{C}(x\otimes1)$
for all $x\in H^{*}(BSO_{4};\mathbb{Z})$ and $\tilde{\beta}_{C}$
is the same as $\beta$ with the only different being the codomain.
With a slight abuse of notation with $\mu$, $\beta$ and $\tilde{\beta}$
being reused, the motivic version of these constructions can be summarized
as the following (\cite{MVW06} pg. 27):

$$\xymatrixcolsep{3.5pc}\xymatrix{
\ & H^{*,*'}(BSO_{4};\mathbb{Z}) \ar[d]_{\otimes 1} \ar[dr]^{\mu} & \ & H^{*+1,*'}(BSO_{4};\mathbb{Z}) \\
0 \ar[r] & H^{*,*'}(BSO_{4};\mathbb{Z})\otimes \mathbb{Z}_{2} \ar[r]^-{\mu_{M}} & H^{*,*'}(BSO_{4};\mathbb{Z}_{2}) \ar[r]^-{\tilde{\beta}_{M}} \ar[ur]^{\tilde{\beta}} \ar[dr]^{\beta} & _{2}H^{*+1,*'}(BSO_{4};\mathbb{Z}) \ar@{^{(}->}[u] \ar[d]^{\mu} \ar[r] & 0 \\
\ & \ & \ & H^{*+1,*'}(BSO_{4};\mathbb{Z}_{2})
}$$

\noindent Note that the subscript $C$'s and $M$'s are meant to distinguish
between ``classical'' and ``motivic''. All together, these sequences
form the following commutative diagram:

$$\xymatrixcolsep{3.5pc}\xymatrix{
0 \ar[r] & H^{*,*'}(BSO_{4};\mathbb{Z})\otimes \mathbb{Z}_{2} \ar[r]^-{\mu_{M}} \ar[d]^{t_{1}} \ar@{-->}[dr] & H^{*,*'}(BSO_{4};\mathbb{Z}_{2}) \ar[r]^-{\tilde{\beta}_{M}} \ar[d]^{t_{2}} \ar@{-->}[dr] & _{2}H^{*+1,*'}(BSO_{4};\mathbb{Z}) \ar[r] \ar[d]^{t_{3}} & 0 \\
0 \ar[r] & H^{*}(BSO_{4};\mathbb{Z})\otimes \mathbb{Z}_{2} \ar[r]^-{\mu_{C}} & H^{*}(BSO_{4};\mathbb{Z}_{2}) \ar[r]^-{\tilde{\beta}_{C}} & _{2}H^{*+1}(BSO_{4};\mathbb{Z}) \ar[r] & 0
}$$

\noindent The maps $t_{i}$ are variants of the realization map $t^{m,n}:H^{m,n}(BSO_{4};\mathbb{Z})\rightarrow H^{m}(BSO_{4};\mathbb{Z})$
where $t=\oplus_{m,n}t^{m,n}$.

There is an important conclusion to draw here that is relevant in
the proof of the main theorem of this section (see \ref{thm: BSOn integral}
below). The map $\mu_{C}$ is injective, so $\varphi\in Ker(\mu_{C}\circ t_{1})$
iff $\varphi\in Ker(t_{1})$. Similarly, $\mu_{M}$ is injective so
$\varphi\in Ker(t_{2}\circ\mu_{M})$ iff $\mu_{M}(\varphi)\in Ker(t_{2})$.
But $\mu_{C}\circ t_{1}=t_{2}\circ\mu_{M}$ by commutativity and therefore

\begin{equation}
\varphi\in Ker(t_{1})\iff\mu_{M}(\varphi)\in Ker(t_{2})
\end{equation}

\begin{thm}
\label{thm: BSOn integral}Let the $d_{i}$'s be Chern classes of
the standard representation of $SO_{4}$ (\cite{Gui07}). Let the
$w_{2i}$'s be the even Stiefel-Whitney classes of the inclusion $G_{2}\subset SO_{7}$
where $deg(w_{2i})=(2i,i+1)$ (\cite{Yag10}). Let the class $y_{2}$
be a lift of twice the classical Euler class of universal $4$-plane
bundle over $BSO_{4}$ (\cite{Fie12}). Take $\tilde{\beta}:H^{*,*'}(BSO_{4};\mathbb{Z}_{2})\rightarrow H^{*+1,*'}(BSO_{4};\mathbb{Z})$
to be the integral Bockstein homomorphism. Finally, let $\tau$ be
a generator of $H^{0,1}(Spec(\mathbb{C});\mathbb{Z}_{2})\text{\ensuremath{\cong\mathbb{Z}_{2}}}$.
Recalling that $H^{*}(BSO_{4};\mathbb{Z})=\mathbb{Z}[p_{1},\sqrt{p_{2}},\tilde{\beta}w_{2}]/(2\tilde{\beta}w_{2})$
for $SO_{4}$ defined over $\mathbb{C}$ (\cite{Bro82}), the integral
motivic cohomology ring of $BSO_{4}$ is

\[
H^{*,*'}(BSO_{4};\mathbb{Z})=H^{*,*'}(Spec(\mathbb{C});\mathbb{Z})[y_{2},d_{2},d_{3},d_{4},\tilde{\beta}\tau^{k}w_{2},\tilde{\beta}\tau^{k-1}w_{2}w_{4}:k\geq0]/I
\]

\noindent where $I$ is the ideal generated by the following relations:

\[
\begin{array}{c}
\begin{array}{ccc}
2d_{3}=0 & y_{2}d_{3}=0 & y_{2}^{2}-4d_{4}=0\end{array}\\
\begin{array}{cccc}
2\tilde{\beta}\tau^{k}w_{2}=0 & 2\tilde{\beta}\tau^{k-1}w_{2}w_{4}=0 & y_{2}\tilde{\beta}\tau^{k}w_{2}=0 & y_{2}\tilde{\beta}\tau^{k-1}w_{2}w_{4}=0\end{array}\\
\begin{array}{c}
(\tilde{\beta}\tau^{k_{1}}w_{2})(\tilde{\beta}\tau^{k_{2}}w_{2})=(\tilde{\beta}\tau^{k_{3}}w_{2})(\tilde{\beta}\tau^{k_{4}}w_{2})\iff k_{1}+k_{2}=k_{3}+k_{4}\\
(\tilde{\beta}\tau^{k_{1}-1}w_{2}w_{4})(\tilde{\beta}\tau^{k_{2}-1}w_{2}w_{4})=(\tilde{\beta}\tau^{k_{3}-1}w_{2}w_{4})(\tilde{\beta}\tau^{k_{4}-1}w_{2}w_{4})\iff k_{1}+k_{2}=k_{3}+k_{4}\\
(\tilde{\beta}\tau^{k_{1}}w_{2})(\tilde{\beta}\tau^{k_{2}-1}w_{2}w_{4})=(\tilde{\beta}\tau^{k_{3}}w_{2})(\tilde{\beta}\tau^{k_{4}-1}w_{2}w_{4})\iff k_{1}+k_{2}=k_{3}+k_{4}\\
\begin{array}{cc}
(\tilde{\beta}\tau^{k-1}w_{2}w_{4})^{2}=d_{4}(\tilde{\beta}\tau^{k}w_{2})^{2} & (\tilde{\beta}\tau^{k}w_{2})^{3}=d_{3}\tilde{\beta}\tau^{3k+1}w_{2}\end{array}
\end{array}
\end{array}
\]

\noindent The degrees of the generators are as follows:

\[
\begin{array}{c}
\begin{array}{cc}
deg(d_{i})=(2i,i) & deg(y_{2})=(4,2)\end{array}\\
\begin{array}{cc}
deg(\tilde{\beta}\tau^{k}w_{2})=(3,2+k) & deg(\tilde{\beta}\tau^{k-1}w_{2}w_{4})=(7,4+k)\end{array}
\end{array}
\]

\noindent Furthermore, the images of the generators under the realization
map $t:H^{*,*'}(BSO_{4};\mathbb{Z})\rightarrow H^{*}(BSO_{4};\mathbb{Z})$
are $t(d_{2})=-p_{1}$, $t(d_{3})=(\tilde{\beta}w_{2})^{2}$, $t(d_{4})=p_{2}$,
$t(y_{2})=2\sqrt{p_{2}}$, $t(\tilde{\beta}\tau^{k}w_{2})=\tilde{\beta}w_{2}$,
and $t(\tilde{\beta}\tau^{k-1}w_{2}w_{4})=\sqrt{p_{2}}\tilde{\beta}w_{2}$.
\end{thm}

\noindent \textbf{Proof:} Note that the integral lifts are sometimes
in bold font for emphasis in diagrams.

The below consideration of only $2$-torsion is sufficient by a transfer
argument. We know that $SL_{2}\times SL_{2}$ is a $2$-fold cover
of $SO_{4}$ and that $H^{*,*'}(BSL_{2};\mathbb{Z})=H^{*,*'}(Spec(\mathbb{C});\mathbb{Z})[c_{2}]$
does not contain torsion outside of $H^{*,*'}(Spec(\mathbb{C});\mathbb{Z})$.
Therefore, the only torsion that can appear in the generators of $H^{*,*'}(BSO_{4};\mathbb{Z})$
over $H^{*,*'}(Spec(\mathbb{C});\mathbb{Z})$ is $2$-torsion.

This method relies on the fact that $H^{*}(BSO_{4};\mathbb{Z})$,
$H^{*}(BSO_{4};\mathbb{Z}_{2})$, $H^{*,*'}(BSO_{4};\mathbb{Z}_{2})$,
and the realization map $t_{2}:H^{*,*'}(BSO_{4};\mathbb{Z}_{2})\rightarrow H^{*}(BSO_{4};\mathbb{Z}_{2})$.
Together with the universal coefficient theorem and $CH^{*}(BSO_{4})$,
these allow us to compute $H^{*,*'}(BSO_{4};\mathbb{Z})$ in its entirety.

$$\xymatrixcolsep{3.5pc}\xymatrix{
0 \ar[r] & H^{*,*'}(BSO_{4};\mathbb{Z})\otimes \mathbb{Z}_{2} \ar[r]^-{\mu_{M}} \ar[d]^{t_{1}} \ar@{-->}[dr] & H^{*,*'}(BSO_{4};\mathbb{Z}_{2}) \ar[r]^-{\tilde{\beta}_{M}} \ar[d]^{t_{2}} \ar@{-->}[dr] & _{2}H^{*+1,*'}(BSO_{4};\mathbb{Z}) \ar[r] \ar[d]^{t_{3}} & 0 \\
0 \ar[r] & H^{*}(BSO_{4};\mathbb{Z})\otimes \mathbb{Z}_{2} \ar[r]^-{\mu_{C}} & H^{*}(BSO_{4};\mathbb{Z}_{2}) \ar[r]^-{\tilde{\beta}_{C}} & _{2}H^{*+1}(BSO_{4};\mathbb{Z}) \ar[r] & 0
}$$

\noindent To fill in the lower left entry, tensoring $H^{*}(BSO_{4};\mathbb{Z})$
with $\mathbb{Z}_{2}$ gives $H^{*}(BSO_{4};\mathbb{Z})\otimes\mathbb{Z}_{2}\cong\mathbb{Z}_{2}[p_{1},\sqrt{p_{2}},\tilde{\beta}w_{2}]$.
\begin{lem}
\label{lem: classical BSO4 exact sequence}In the exact sequence $0\rightarrow H^{*}(BSO_{4};\mathbb{Z})\otimes\mathbb{Z}_{2}\xrightarrow{\mu_{C}}H^{*}(BSO_{4};\mathbb{Z}_{2})\xrightarrow{\tilde{\beta}_{C}}{}_{2}H^{*+1}(BSO_{4};\mathbb{Z})\rightarrow0$:
\end{lem}

\[
\begin{array}{c}
\begin{array}{cc}
\mu_{C}(p_{1}^{k_{1}}(\sqrt{p_{2}})^{k_{2}}(\tilde{\beta}w_{2})^{k_{3}}\otimes1)=w_{2}^{2k_{1}}w_{3}^{k_{3}}w_{4}^{k_{2}} & \tilde{\beta}_{C}(w_{2}^{2k_{1}}w_{3}^{k_{3}}w_{4}^{k_{2}})=0\end{array}\\
\tilde{\beta}_{C}(w_{2}^{2k_{1}+1}w_{3}^{k_{3}}w_{4}^{k_{2}})=p_{1}^{k_{1}}(\sqrt{p_{2}})^{k_{2}}(\tilde{\beta}w_{2})^{k_{3}+1}
\end{array}
\]

\noindent \textbf{Proof of Lemma:} Brown's Theorem $1.5$ includes
the statement that $\mu(p_{i})=w_{2i}^{2}$ for all $i$. This shows
that $\mu_{C}(p_{1}\otimes1)=w_{2}^{2}$ and $\mu_{C}(\sqrt{p_{2}}\otimes1)=w_{4}$.
Injectivity of $\mu_{C}$ along with the facts that $H^{3}(BSO_{4};\mathbb{Z})=\mathbb{Z}_{2}\{\tilde{\beta}w_{2}\}$
and $H^{3}(BSO_{4};\mathbb{Z}_{2})=\mathbb{Z}_{2}\{w_{3}\}$ gives
that $\mu_{C}(\tilde{\beta}w_{2}\otimes1)=w_{3}$. The first line
follows from exactness.

The approach for the second line is different because $\tilde{\beta}_{C}$
is a module homomorphism, not a ring homomorphism. Clearly $\tilde{\beta}_{C}(w_{2})=\tilde{\beta}w_{2}$,
as indicated by Brown's notation. More generally, we know that $Sq^{1}=\beta=\mu\circ\tilde{\beta}_{C}$
is simply the Bockstein homomorphism and the first Steenrod squaring
operation. On Stiefel-Whitney classes, $\beta$ is characterized by
the fact that $\beta(w_{2i})=w_{2i+1}$ and $\beta(w_{2i+1})=0$ (\cite{HN18});
this along with induction on the Cartan formula and the fact that
$w_{5}=0$ gives the result:

\[
\begin{array}{c}
\beta(w_{2}^{2k_{1}+1})=w_{2}\beta(w_{2}^{2k_{1}})+\beta(w_{2})w_{2}^{2k_{1}}=w_{2}^{2k_{1}}w_{3}=\mu(p_{1}^{k_{1}}\tilde{\beta}w_{2})\\
\beta(w_{2}^{2k_{1}+1}w_{3}^{k_{3}})=w_{2}^{2k_{1}+1}\beta(w_{3}^{k_{3}})+\beta(w_{2}^{2k_{1}+1})w_{3}^{k_{3}}=w_{2}^{2k_{1}}w_{3}^{k_{3}+1}=\mu(p_{1}^{k_{1}}(\tilde{\beta}w_{2})^{k_{3}+1})\\
\beta(w_{2}^{2k_{1}+1}w_{3}^{k_{3}}w_{4}^{k_{2}})=w_{2}^{2k_{1}+1}w_{3}^{k_{3}}\beta(w_{4}^{k_{2}})+\beta(w_{2}^{2k_{1}+1}w_{3}^{k_{3}})w_{4}^{k_{2}}=w_{2}^{2k_{1}}w_{3}^{k_{3}+1}w_{4}^{k_{2}}=\mu(p_{1}^{k_{1}}(\sqrt{p_{2}})^{k_{2}}(\tilde{\beta}w_{2})^{k_{3}+1})
\end{array}
\]

\noindent \begin{flushright}
$\boxempty$
\par\end{flushright}

The next step is to consider the elements of $CH^{*}(BSO_{4})$ as
lifts of classes in $H^{*}(BSO_{4};\mathbb{Z})$; in particular, we
need to contextualize $d_{2}$, $d_{3}$, $d_{4}$ and $y_{2}$. It
well-known that Pontryagin classes have lifts to the motivic setting
via $t(d_{2i})=(-1)^{i}p_{i}$ (\cite{Bro82}), so the following diagrams
are known from \cite{Yag10}:

$$\xymatrix@!C=8pc{
\bm{d_{2}} \ar@{|->}[r]^-{\mu} \ar@{|->}[d]^{t} & \tau^{-2}(w_{2}^{alg})^{2} \ar@{|->}[d]^{t_{2}} & \bm{d_{3}} \ar@{|->}[r]^-{\mu} \ar@{|->}[d]^{t} & \tau^{-1}(w_{3}^{alg})^{2}  \ar@{|->}[d]^{t_{2}} \\
-p_{1} \ar@{|->}[r]^-{\mu} & w_{2}^{2} &  (\tilde{\beta}w_{2})^{2} \ar@{|->}[r]^-{\mu} & w_{3}^{2} \\
\bm{d_{4}} \ar@{|->}[r]^-{\mu} \ar@{|->}[d]^{t} & \tau^{-2}(w_{4}^{alg})^{2} \ar@{|->}[d]^{t_{2}} & \bm{y_{2}} \ar@{|->}[r]^-{\mu} \ar@{|->}[d]^{t} & y_{0,2}  \ar@{|->}[d]^{t_{2}} \\
p_{2} \ar@{|->}[r]^-{\mu} & w_{4}^{2} &  2\sqrt{p_{2}} \ar@{|->}[r]^-{\mu} & 0
}$$

\noindent On the other hand, there is no lift $e^{alg}$ of the Euler
class $e=\sqrt{p_{2}}$. We know from Field's Theorem 1 that her class
$y_{2}\in CH^{*}(BSO_{4})$ is mapped to $2\sqrt{p_{2}}$ under the
realization map; clearly this implies that twice the Euler class does
have a lift in $H^{*,*'}(BSO_{4};\mathbb{Z})$. We run into trouble
if we instead attempt to lift only $\sqrt{p_{2}}$, as seen in detail
below:
\begin{lem}
\label{lem: no lift of Euler class}There is no lift of $\sqrt{p_{2}}\in H^{*}(BSO_{4};\mathbb{Z})$
to $H^{*,*'}(BSO_{4};\mathbb{Z})$.
\end{lem}

\noindent \textbf{Proof of Lemma:} Suppose such of lift $e^{alg}\in H^{4,n}(BSO_{4};\mathbb{Z})$
of $\sqrt{p_{2}}$ exists for some $n\geq2$.
\begin{itemize}
\item Clearly $(\sqrt{p_{2}})^{2}-p_{2}=0$; thus $((e^{alg})^{2}-d_{4})\otimes1\in Ker(t_{1})$
and so $\mu_{M}(((e^{alg})^{2}-d_{4})\otimes1)\in Ker(t_{2})$.
\item Simplifying gives $\mu_{M}(((e^{alg})^{2}-d_{4})\otimes1)=\mu_{M}(e^{alg}\otimes1)^{2}-\mu(d_{4})$.
\item We know from Yagita's work that $H^{8,4}(BSO_{4};\mathbb{Z}_{2})=\mathbb{Z}_{2}\{\mu(d_{2})^{2},y_{0,2}\mu(d_{2}),\mu(d_{4})\}$
and if $k\geq5$ then $H^{8,k}(BSO_{4};\mathbb{Z}_{2})=\mathbb{Z}_{2}\{\tau^{k-5}w_{4}^{alg}\mu(d_{2}),\tau^{n-5}w_{2}^{alg}\mu(d_{3}),\tau^{n-4}\mu(d_{2})^{2},\tau^{n-4}\mu(d_{4})\}$.
\item Computations show that $t_{2}$ sends the following elements of $H^{*,*'}(BSO_{4};\mathbb{Z}_{2})$
to their counterparts in $H^{*}(BSO_{4};\mathbb{Z}_{2}$): $\tau\mapsto1$,
$\mu(d_{2})^{2}\mapsto w_{2}^{4}$, $y_{0,2}\mu(d_{2})\mapsto0$,
$\mu(d_{4})\mapsto w_{4}^{2}$, $w_{4}^{alg}\mu(d_{2})\mapsto w_{2}^{2}w_{4}$,
and $w_{2}^{alg}\mu(d_{3})\mapsto w_{2}w_{3}^{2}$. Thus, $\mu_{M}(e^{alg}\otimes1)^{2}\in H^{8,2n}(BSO_{4};\mathbb{Z}_{2})$
must be equal to $\tau^{2n-4}\mu(d_{4})$ so that $t_{2}(\mu_{M}(e^{alg}\otimes1)^{2})=w_{4}^{2}$.
\item This implies that if $e^{alg}$ exists then it must satisfy the property
that $\mu_{M}(e^{alg}\otimes1)$ is equal to $\tau^{n_{3}-2}x$ for
some $x\in H^{4,2}(BSO_{4};\mathbb{Z}_{2})$ such that $x^{2}=\mu(d_{4})$.
The problem is that no such element $x$ exists because $H^{4,2}(BSO_{4};\mathbb{Z}_{2})=\mathbb{Z}_{2}\{\mu(d_{2}),y_{0,2}\}$
and neither of these two generators satisfy this property.
\item Thus, the lift $e^{alg}$ cannot exist.
\end{itemize}
Visually, the above argument can be summarized as follows:

$$\xymatrix@!C=11pc{
\ & w_{4}^{alg} \ar@{|->}[r]^-{\tilde{\beta}_{M}} \ar@{|->}[d]^{t_{2}} & 0 \ar@{|->}[d]^{t} \\
\sqrt{p_{2}} \ar@{|->}[r]^-{\mu} & w_{4} \ar@{|->}[r]^-{\tilde{\beta}_{C}} & 0
}$$
\noindent \begin{flushright}
$\boxempty$
\par\end{flushright}

The key to the remainder of this proof is that whenever $_{2}H^{m+1}(BSO_{4};\mathbb{Z})$
is non-zero there will be elements of $H^{m+1}(BSO_{4};\mathbb{Z})$
that can be lifted to non-zero elements of $H^{m+1,n}(BSO_{4};\mathbb{Z})$
for some $n$. To see this, note that we can rewrite $H^{*}(BSO_{4};\mathbb{Z})$
in the following way:

\[
H^{*}(BSO_{4};\mathbb{Z})=\mathbb{Z}\{p_{1}^{i}(\sqrt{p_{2}})^{j}:i,j\in\mathbb{N}\}\oplus\mathbb{Z}_{2}\{p_{1}^{i}(\sqrt{p_{2}})^{j}(\tilde{\beta}w_{2})^{k}:i,j\in\mathbb{N},k\in\mathbb{Z}_{>0}\}
\]

\noindent So $_{2}H^{m+1}(BSO_{4};\mathbb{Z})=\mathbb{Z}_{2}\{p_{1}^{i}(\sqrt{p_{2}})^{j}(\tilde{\beta}w_{2})^{k}:i,j\in\mathbb{N},k\in\mathbb{Z}_{>0},4i+4j+3k=m+1\}$
is non-trivial precisely when $[\frac{1}{4}(m+2),\frac{1}{3}(m+1)]\cap\mathbb{Z}\neq\emptyset$;
a quick calculation shows that this is equivalent to having $m+1\in\{3,6,7\}\cup\mathbb{Z}_{\geq9}$.
Conversely, we now know that $_{2}H^{m+1}(BSO_{4};\mathbb{Z})=\{0\}$
only if $m+1\in\{1,2,4,5,8\}$.

Consider the values of $0\rightarrow H^{m}(BSO_{4};\mathbb{Z})\otimes\mathbb{Z}_{2}\xrightarrow{\mu_{C}}H^{m}(BSO_{4};\mathbb{Z}_{2})\xrightarrow{\tilde{\beta}_{C}}_{2}H^{m+1}(BSO_{4};\mathbb{Z})\rightarrow0$
for the first several values of $m$, as given by the Lemma \ref{lem: classical BSO4 exact sequence}:

$$\xymatrixcolsep{2.3pc}\xymatrixrowsep{0.6pc}\xymatrix{
m=0: & 0 \ar[r] & \mathbb{Z}_{2}\{1\otimes1\} \ar[r] & \mathbb{Z}_{2}\{1\} \ar[r] & 0 \ar[r] & 0 \\
m=1: & 0 \ar[r] & 0 \ar[r] & 0 \ar[r] & 0 \ar[r] & 0 \\
m=2: & 0 \ar[r] & 0 \ar[r] & \mathbb{Z}_{2}\{w_{2}\} \ar[r] & \mathbb{Z}_{2}\{\tilde{\beta}w_{2}\} \ar[r] & 0 \\
m=3: & 0 \ar[r] & \mathbb{Z}_{2}\{\tilde{\beta}w_{2}\otimes1\} \ar[r] & \mathbb{Z}_{2}\{w_{3}\} \ar[r] & 0 \ar[r] & 0 \\
m=4: & 0 \ar[r] & \mathbb{Z}_{2}\{p_{1}\otimes1,\sqrt{p_{2}}\otimes1\} \ar[r] & \mathbb{Z}_{2}\{w_{2}^{2},w_{4}\} \ar[r] & 0 \ar[r] & 0 \\
m=5: & 0 \ar[r] & 0 \ar[r] & \mathbb{Z}_{2}\{w_{2}w_{3}\} \ar[r] & \mathbb{Z}_{2}\{(\tilde{\beta}w_{2})^{2}\} \ar[r] & 0 \\
m=6: & 0 \ar[r] & \mathbb{Z}_{2}\{(\tilde{\beta}w_{2})^{2}\otimes1\} \ar[r] & \mathbb{Z}_{2}\{w_{2}^{3},w_{2}w_{4},w_{3}^{2}\} \ar[r] & \mathbb{Z}_{2}\{p_{1}\tilde{\beta}w_{2},\sqrt{p_{2}}\tilde{\beta}w_{2}\} \ar[r] & 0 \\
}$$

\noindent Keeping in mind that $\tilde{\beta}_{C}\circ t_{2}=t_{3}\circ\tilde{\beta}_{M}$:
\begin{itemize}
\item The $m=2$ case shows that $\tilde{\beta}_{C}(w_{2})=\tilde{\beta}w_{2}$.
Because $t_{2}(\tau^{k}w_{2}^{alg})=w_{2}$ for any choice of $k$,
we have that $t_{3}\circ\tilde{\beta}_{M}(\tau^{k}w_{2}^{alg})=\tilde{\beta}w_{2}$.
This means that $\tilde{\beta}_{M}(\tau^{k}w_{2}^{alg})$ is non-zero
in $_{2}H^{3,2+i}(BSO_{4};\mathbb{Z})$, i.e. $\tilde{\beta}_{M}(\tau^{i}w_{2}^{alg})\neq0$
but $2\tilde{\beta}_{M}(\tau^{i}w_{2}^{alg})=0$. In other words,
$\tilde{\beta}_{M}(\tau^{k}w_{2}^{alg})$ is a lift of $\tilde{\beta}w_{2}$
for all $k$.
\begin{itemize}
\item Note that this demonstrates the distinction between the algebraic
lifts $w_{i}^{alg}$ and the topological lifts $w_{i}^{top}$ in $H^{*,*'}(BSO_{4};\mathbb{Z}_{2})$.
Both give rise to lifts in $H^{*,*'}(BSO_{4};\mathbb{Z})$ and, in
fact, so do infinitely many other classes:
\end{itemize}
$$\xymatrixcolsep{3pc}\xymatrix{
w_{2}^{alg} \ar@{|->}[r]^-{\tilde{\beta}_{M}} \ar@{|->}[d]^{t_{2}} & \tilde{\beta}_{M}(w_{2}^{alg}) \ar@{|->}[r]^-{\mu_{M}} \ar@{|->}[d]^{t_{3}} & w_{3}^{alg} \ar@{|->}[d]^{t_{2}} &
\tau w_{2}^{alg} \ar@{|->}[r]^-{\tilde{\beta}_{M}} \ar@{|->}[d]^{t_{2}} & \tilde{\beta}_{M}(\tau w_{2}^{alg}) \ar@{|->}[r]^-{\mu_{M}} \ar@{|->}[d]^{t_{3}} & w_{3}^{top} \ar@{|->}[d]^{t_{2}} \\
w_{2} \ar@{|->}[r]^-{\tilde{\beta}_{C}} & \tilde{\beta}w_{2} \ar@{|->}[r]^-{\mu_{C}} & w_{3} &
w_{2} \ar@{|->}[r]^-{\tilde{\beta}_{C}} & \tilde{\beta}w_{2} \ar@{|->}[r]^-{\mu_{C}} & w_{3}
}$$
\item Similarly, the $m=5$ case shows that $\tilde{\beta}_{M}(\tau^{k}w_{2}^{alg}w_{3}^{alg})$
is a lift of $(\tilde{\beta}w_{2})^{2}$ for all $k$.
\item On the other hand, the Cartan formula gives that $\beta(\tau^{k}w_{3}^{alg})=0$;
it is well-known that $\beta(\tau^{k})=0$ because no integral analog
of $\tau$ exists (\cite{MVW06} 4.2) and we also know that $\beta(w_{3}^{alg})=0$
(\cite{HN18}). This shows that $\beta(\tau^{k}w_{3}^{alg})$ is not
a non-trivial lift of any element of $H^{4}(BSO_{4};\mathbb{Z})$.
\item The same reasoning applied to the $m=4$ case shows that $\tilde{\beta}_{M}(\tau^{k}w_{2}^{2})$
and $\tilde{\beta}_{M}(\tau^{k}w_{4})$ are both $0$ and thus are
not considered lifts either.
\end{itemize}
More generally, all $2$-torsion elements in $H^{*}(BSO_{4};\mathbb{Z})$
will have lifts to $2$-torsion elements in $H^{*,*'}(BSO_{4};\mathbb{Z})$.
As above, these lifts will be of the form $\tilde{\beta}(\tau^{k}x)$
for classes $x\in H^{*,*'}(BSO_{4};\mathbb{Z}_{2})$ such that $\tilde{\beta}_{C}\circ t_{2}(x)$
is $2$-torsion.

Ultimately the result of the above argument is that the integral motivic
cohomology of $BSO_{4}$ has infinitely many generators. For example,
each $\tilde{\beta}_{M}(\tau^{k}w_{2}^{alg})\in H^{3,2+k}(BSO_{4};\mathbb{Z})$
is $2$-torsion lift satisfying the property $\tau^{k}w_{3}^{alg}=\mu\circ\tilde{\beta}_{M}(\tau^{k}w_{2}^{alg})$.
It is key to determine the relations between these $\tilde{\beta}_{M}(\tau^{k}w_{2}^{alg})$
for different $k$. One immediate relation mod $2$ is given by

\begin{equation}
k\geq k'\implies\mu\circ\tilde{\beta}_{M}(\tau^{k}w_{2}^{alg})=\tau^{k-k'}\mu\circ\tilde{\beta}_{M}(\tau^{k'}w_{2}^{alg})
\end{equation}

\noindent It would be convenient if there were a similar relation
on the integral side of things, i.e. $k\geq k'$ implies $\tilde{\beta}_{M}(\tau^{k}w_{2})=T^{k-k'}\tilde{\beta}_{M}(\tau^{k'}w_{2})$
for some $T\in H^{0,1}(BSO_{4};\mathbb{Z})$. Unfortunately that class
would have to live in $H^{0,1}(Spec(\mathbb{C});\mathbb{Z})$ and
it is known that no such class exists (\cite{MVW06} 4.2).

Clearly the integral case is more complicated than the mod $2$ case
due to the fact that no integral analog of $\tau$ exists. Before
proceeding any further, the following lemma will help to greatly simplify
arguments moving forward:
\begin{lem}
\label{lem: mod 2 reductions equal means originals equal}Generators
of $2$-torsion in $H^{*,*'}(BSO_{4};\mathbb{Z})$ that are not $2$-divisible
are equal iff their mod $2$ reductions in $H^{*,*'}(BSO_{4};\mathbb{Z}_{2})$
are equal.
\end{lem}

\noindent \textbf{Proof of Lemma:} Suppose $x,y\in H^{*,*'}(BSO_{4};\mathbb{Z})$
are $2$-torsion and are not $2$-divisible. Clearly if $x=y$ then
$\mu(x)=\mu(y)$, so it remains to go in the other direction. The
elements $x$ and $y$ are not $2$-divisible so $x\otimes1$ and
$y\otimes1$ are both non-zero in $H^{*,*'}(BSO_{4};\mathbb{Z})\otimes\mathbb{Z}_{2}$. 

Assume that $\mu(x)=\mu(y)$. This implies that $\mu_{M}(x\otimes1)=\mu_{M}(y\otimes1)$.
Injectivity of $\mu_{M}$ gives that $x\otimes1=y\otimes1$. Thus,
some odd multiples of $x$ and $y$ must be equal, i.e. $(2i+1)x=(2j+1)y$
for some $i$ and $j$. But $x$ and $y$ are $2$-torsion and therefore
$2ix$ and $2jy$ are zero; we conclude that $x=y$.
\noindent \begin{flushright}
$\boxempty$
\par\end{flushright}

We have shown that there are infinitely many lifts of $\tilde{\beta}w_{2}$;
as stated above, in principal any $2$-torsion element of $H^{*}(BSO_{4};\mathbb{Z})$
will have a lift in a similar form. Applying the same construction
as done for $\tilde{\beta}w_{2}$ to other elements of Yagita's presentation
yields all remaining lifts (where here we make a slight simplification
to the notation of the lifts):
\begin{flushleft}
$$\xymatrix@!C=11pc{
\tau^{k}w_{2}^{alg} \ar@{|->}[r]^-{\tilde{\beta}_{M}} \ar@{|->}[d]^{t_{2}} & \bm{\tilde{\beta}\tau^{k}w_{2}} \ar@{|->}[r]^-{\mu} \ar@{|->}[d]^{t} & \tau^{k}w_{3}^{alg} \ar@{|->}[d]^{t_{2}} \\
w_{2} \ar@{|->}[r]^-{\tilde{\beta}_{C}} & \tilde{\beta}w_{2} \ar@{|->}[r]^-{\mu} & w_{3}
}$$
\par\end{flushleft}

\begin{flushleft}
$$\xymatrix@!C=11pc{
\tau^{k-1}w_{2}^{alg}w_{3}^{alg} \ar@{|->}[r]^-{\tilde{\beta}_{M}} \ar@{|->}[d]^{t_{2}} & \bm{\tilde{\beta}\tau^{k-1}w_{2}w_{3}} \ar@{|->}[r]^-{\mu} \ar@{|->}[d]^{t} & \tau^{k-1}(w_{3}^{alg})^{2} \ar@{|->}[d]^{t_{2}} \\
w_{2}w_{3} \ar@{|->}[r]^-{\tilde{\beta}_{C}} & (\tilde{\beta}w_{2})^{2} \ar@{|->}[r]^-{\mu} & w_{3}^{2}
}$$
\par\end{flushleft}

\begin{flushleft}
$$\xymatrix@!C=11pc{
\tau^{k-1}w_{2}^{alg}w_{4}^{alg} \ar@{|->}[r]^-{\tilde{\beta}_{M}} \ar@{|->}[d]^{t_{2}} & \bm{\tilde{\beta}\tau^{k-1}w_{2}w_{4}} \ar@{|->}[r]^-{\mu} \ar@{|->}[d]^{t} & \tau^{k-1}w_{3}^{alg}w_{4}^{alg} \ar@{|->}[d]^{t_{2}} \\
w_{2}w_{4} \ar@{|->}[r]^-{\tilde{\beta}_{C}} & \sqrt{p_{2}}\tilde{\beta}w_{2} \ar@{|->}[r]^-{\mu} & w_{3}w_{4}
}$$
\par\end{flushleft}

\begin{flushleft}
$$\xymatrix@!C=11pc{
\tau^{k-1}w_{2}^{alg}w_{3}^{alg}w_{4}^{alg} \ar@{|->}[r]^-{\tilde{\beta}_{M}} \ar@{|->}[d]^{t_{2}} & \bm{\tilde{\beta}\tau^{k-1}w_{2}w_{3}w_{4}} \ar@{|->}[r]^-{\mu} \ar@{|->}[d]^{t} & \tau^{k-1}(w_{3}^{alg})^{2}w_{4}^{alg} \ar@{|->}[d]^{t_{2}} \\
w_{2}w_{3}w_{4} \ar@{|->}[r]^-{\tilde{\beta}_{C}} & \sqrt{p_{2}}(\tilde{\beta}w_{2})^{2} \ar@{|->}[r]^-{\mu} & w_{3}^{2}w_{4}
}$$
\par\end{flushleft}

\begin{flushleft}
$$\xymatrix@!C=11pc{
\tau^{k-2}w_{2}^{alg}(w_{4}^{alg})^{2} \ar@{|->}[r]^-{\tilde{\beta}_{M}} \ar@{|->}[d]^{t_{2}} & \bm{\tilde{\beta}\tau^{k-2}w_{2}w_{4}^{2}} \ar@{|->}[r]^-{\mu} \ar@{|->}[d]^{t} & \tau^{k-2}w_{3}^{alg}(w_{4}^{alg})^{2} \ar@{|->}[d]^{t_{2}} \\
w_{2}w_{4}^{2} \ar@{|->}[r]^-{\tilde{\beta}_{C}} & p_{2}\tilde{\beta}w_{2} \ar@{|->}[r]^-{\mu} & w_{3}w_{4}^{2}
}$$
\par\end{flushleft}

\begin{flushleft}
$$\xymatrix@!C=11pc{
\tau^{k-3}w_{2}^{alg}w_{3}^{alg}(w_{4}^{alg})^{2} \ar@{|->}[r]^-{\tilde{\beta}_{M}} \ar@{|->}[d]^{t_{2}} & \bm{\tilde{\beta}\tau^{k-3}w_{2}w_{3}w_{4}^{2}} \ar@{|->}[r]^-{\mu} \ar@{|->}[d]^{t} & \tau^{k-3}(w_{3}^{alg})^{2}(w_{4}^{alg})^{2} \ar@{|->}[d]^{t_{2}} \\
w_{2}w_{3}w_{4}^{2} \ar@{|->}[r]^-{\tilde{\beta}_{C}} & p_{2}(\tilde{\beta}w_{2})^{2} \ar@{|->}[r]^-{\mu} & w_{3}^{2}w_{4}^{2}
}$$
\par\end{flushleft}

\noindent One important conclusion to draw from these diagrams is
that the generator $d_{3}$ is of this form:

\begin{equation}
d_{3}=\tilde{\beta}\tau^{-1}w_{2}w_{3}
\end{equation}

\noindent There are two other types lift that we have seen, namely
those involving $d_{2}$, $d_{4}$, and $y_{2}$. However, these types
do not involve $\tau$ because multiple lifts can only exist for $2$-torsion
elements in this construction:

$$\xymatrix@!C=11pc{
\bm{d_{2}^{k_1}d_{4}^{k_{2}}} \ar@{|->}[r]^-{\mu} \ar@{|->}[d]^{t} & \mu(d_{2})^{k_{1}}\mu(d_{4})^{k_{2}} \ar@{|->}[d]^{t_{2}} \\
(-p_{1})^{k_{1}}(\sqrt{p_{2}})^{2k_{2}} \ar@{|->}[r]^-{\mu} & w_{2}^{2k_{1}}w_{4}^{2k_{2}}
}$$

$$\xymatrix@!C=11pc{
\bm{y_{2}d_{2}^{k_1}d_{4}^{k_2}} \ar@{|->}[r]^-{\mu} \ar@{|->}[d]^{t} & y_{0,2}\mu(d_{2})^{k_{1}}\mu(d_{4})^{k_{2}} \ar@{|->}[d]^{t_{2}} \\
2(-p_{1})^{k_1}(\sqrt{p_{2}})^{2k_2+1} \ar@{|->}[r]^-{\mu} & 0
}$$

If we take ``simple'' to mean a product of only generators with
no coefficients, we see that almost all simple elements of $H^{*}(BSO_{4};\mathbb{Z})$
have lifts. There is only one type of simple element that remains
to be considered, namely the case of the power of $\sqrt{p_{2}}$
being odd and there being no $\tilde{\beta}w_{2}$. The Euler class
falls into this category and so the following lemma is a natural continuation
of the logic from Lemma \ref{lem: no lift of Euler class}:
\begin{lem}
An element $p_{1}^{j_{1}}(\sqrt{p_{2}})^{j_{2}}(\tilde{\beta}w_{2})^{j_{3}}\in H^{*}(BSO_{4};\mathbb{Z})$
lifts $H^{*,*'}(BSO_{4};\mathbb{Z})$ iff $j_{2}$ is even or $j_{3}>0$.
\end{lem}

\noindent \textbf{Proof of Lemma:} We've seen that if $j_{3}$ is
greater than zero then there is a lift. In addition, we know that
$(-p_{1})^{k_{1}}(\sqrt{p_{2}})^{2k_{2}}$ lifts to $d_{2}^{k_{1}}d_{4}^{k_{2}}$.
Thus, it suffices to show that if $j_{2}=2k_{2}+1$ and $k_{3}=0$
then there is no lift:
\begin{itemize}
\item Suppose a lift of $l^{alg}\in H^{4k_{1}+8k_{2}+4,n}(BSO_{4};\mathbb{Z})$
of $p_{1}^{k_{1}}(\sqrt{p_{2}})^{2k_{2}+1}$ exists.
\item Clearly $(p_{1}^{k_{1}}(\sqrt{p_{2}})^{2k_{2}+1})=p_{1}^{2k_{1}}(\sqrt{p_{2}})^{4k_{2}+2}$
because $p_{1}$ and $p_{2}$ are both even degree.
\item But $p_{1}^{2k_{1}}(\sqrt{p_{2}})^{4k_{2}+2}=p_{1}^{2k_{1}}p_{2}^{2k_{2}+1}$
lifts to $d_{2}^{2k_{1}}d_{4}^{2k_{2}+1}$ and so $((l^{alg})^{2}-d_{2}^{2k_{1}}d_{4}^{2k_{2}+1})\otimes1$
is in $Ker(t_{1})$.
\item Thus, $\mu_{M}(((l^{alg})^{2}-d_{2}^{2k_{1}}d_{4}^{2k_{2}+1})\otimes1)=\mu_{M}(l^{alg}\otimes1)^{2}-\mu(d_{2})^{2k_{1}}\mu(d_{4})^{2k_{2}+1}$
is in $Ker(t_{2})$.
\item This gives $\mu_{M}(l^{alg}\otimes1)^{2}=\tau^{2n-4k_{1}-8k_{2}-4}\mu(d_{2})^{2k_{1}}\mu(d_{4})^{2k_{2}+1}$
and there must be some $x\in H^{4,2}(BSO_{4};\mathbb{Z}_{2})$ such
that $x^{2}=\mu(d_{4})$ where $\mu_{M}(l^{alg}\otimes1)$ is equal
to $\tau^{n-2k_{1}-4k_{2}-2}\mu(d_{2})^{k_{1}}\mu(d_{4})^{k_{2}}x$.
\item No such $x$ exists in $H^{4,2}(BSO_{4};\mathbb{Z}_{2})=\mathbb{Z}_{2}\{\mu(d_{2}),y_{0,2}\}$
and so $l^{alg}$ does not exist.
\end{itemize}
\noindent \begin{flushright}
$\boxempty$
\par\end{flushright}

\begin{sidewaystable}
\begin{centering}
\begin{tabular}{|c|c|c|c|c|c|c|}
\hline 
{\small{}$\begin{array}{c}
\\
\\
\end{array}$}Family{\small{}$\begin{array}{c}
\\
\\
\end{array}$} & Form & Coefficient & Power of $p_{1}$ & Power of $\sqrt{p_{2}}$ & Power of $\tilde{\beta}w_{2}$ & Form of Lifts\tabularnewline
\hline 
\hline 
{\small{}$\begin{array}{c}
\\
\\
\end{array}$}$1${\small{}$\begin{array}{c}
\\
\\
\end{array}$} & $p_{1}^{k_{2}}(\tilde{\beta}w_{2})^{2k_{3}+1}$ & $1$ & any & $0$ & odd & $d_{2}^{k_{2}}d_{3}^{k_{3}}\tilde{\beta}\tau^{k_{1}}w_{2}$\tabularnewline
\hline 
{\small{}$\begin{array}{c}
\\
\\
\end{array}$}$2${\small{}$\begin{array}{c}
\\
\\
\end{array}$} & $p_{1}^{k_{2}}(\tilde{\beta}w_{2})^{2k_{3}+2}$ & $1$ & any & $0$ & even $>0$ & $d_{2}^{k_{2}}d_{3}^{k_{3}}\tilde{\beta}\tau^{k_{1}-1}w_{2}w_{3}$\tabularnewline
\hline 
{\small{}$\begin{array}{c}
\\
\\
\end{array}$}$3${\small{}$\begin{array}{c}
\\
\\
\end{array}$} & $p_{1}^{k_{2}}(\sqrt{p_{2}})^{2k_{4}+1}(\tilde{\beta}w_{2})^{2k_{3}+1}$ & $1$ & any & odd & odd & $d_{2}^{k_{2}}d_{3}^{k_{3}}d_{4}^{k_{4}}\tilde{\beta}\tau^{k_{1}-1}w_{2}w_{4}$\tabularnewline
\hline 
{\small{}$\begin{array}{c}
\\
\\
\end{array}$}$4${\small{}$\begin{array}{c}
\\
\\
\end{array}$} & $p_{1}^{k_{2}}(\sqrt{p_{2}})^{2k_{4}+1}(\tilde{\beta}w_{2})^{2k_{3}+2}$ & $1$ & any & odd & even $>0$ & $d_{2}^{k_{2}}d_{3}^{k_{3}}d_{4}^{k_{4}}\tilde{\beta}\tau^{k_{1}-1}w_{2}w_{3}w_{4}$\tabularnewline
\hline 
{\small{}$\begin{array}{c}
\\
\\
\end{array}$}$5${\small{}$\begin{array}{c}
\\
\\
\end{array}$} & $p_{1}^{k_{2}}(\sqrt{p_{2}})^{2k_{4}+2}(\tilde{\beta}w_{2})^{2k_{3}+1}$ & $1$ & any & even $>0$ & odd & $d_{2}^{k_{2}}d_{3}^{k_{3}}d_{4}^{k_{4}}\tilde{\beta}\tau^{k_{1}-2}w_{2}w_{4}^{2}$\tabularnewline
\hline 
{\small{}$\begin{array}{c}
\\
\\
\end{array}$}$6${\small{}$\begin{array}{c}
\\
\\
\end{array}$} & $p_{1}^{k_{2}}(\sqrt{p_{2}})^{2k_{4}+2}(\tilde{\beta}w_{2})^{2k_{3}+2}$ & $1$ & any & even $>0$ & even $>0$ & $d_{2}^{k_{2}}d_{3}^{k_{3}}d_{4}^{k_{4}}\tilde{\beta}\tau^{k_{1}-3}w_{2}w_{3}w_{4}^{2}$\tabularnewline
\hline 
{\small{}$\begin{array}{c}
\\
\\
\end{array}$}$7${\small{}$\begin{array}{c}
\\
\\
\end{array}$} & $\lambda p_{1}^{k_{2}}(\sqrt{p_{2}})^{2k_{4}}$ & any & any & even & $0$ & $\lambda(-d_{2})^{k_{2}}d_{4}^{k_{4}}$\tabularnewline
\hline 
{\small{}$\begin{array}{c}
\\
\\
\end{array}$}$8${\small{}$\begin{array}{c}
\\
\\
\end{array}$} & $2\lambda p_{1}^{k_{2}}(\sqrt{p_{2}})^{2k_{4}+1}$ & even & any & odd & $0$ & $\lambda y_{2}(-d_{2})^{k_{2}}d_{4}^{k_{4}}$\tabularnewline
\hline 
{\small{}$\begin{array}{c}
\\
\\
\end{array}$}$9${\small{}$\begin{array}{c}
\\
\\
\end{array}$} & $(2\lambda+1)p_{1}^{k_{1}}(\sqrt{p_{2}})^{2k_{2}+1}$ & odd & any & odd & $0$ & N/A\tabularnewline
\hline 
\end{tabular}
\par\end{centering}
\caption{\label{tab: table of integral BSO4 lifts}The elements of $H^{*,*'}(BSO_{4};\mathbb{Z})$
fall into nine families. Eight of these have lifts as above but the
ninth does not.}
\end{sidewaystable}

The first objective is to remove any generators that can be rewritten
completely in terms of other generators. All the generators from the
first six families in Table \ref{tab: table of integral BSO4 lifts}
are $2$-torsion but not $2$-divisible and so Lemma \ref{lem: mod 2 reductions equal means originals equal}
applies for finding integral relations. The driving idea for these
relations comes from \cite{Por21} 7.3.1; in particular, if two classes
in $H^{*,*'}(BSO_{4};\mathbb{Z}_{2})-Ker(t_{2})$ have the same realizations
and degrees then they must be equal. This fact allows us to rewrite
families $2$, $4$, $5$ and $6$:

\begin{equation}
\begin{array}{c}
\begin{array}{cc}
\tilde{\beta}\tau^{-1}w_{2}w_{3}=d_{3} & \tilde{\beta}\tau^{k_{1}+k_{2}}w_{2}w_{3}=(\tilde{\beta}\tau^{k_{1}}w_{2})(\tilde{\beta}\tau^{k_{2}}w_{2})\end{array}\\
\tilde{\beta}\tau^{k_{1}+k_{2}-1}w_{2}w_{3}w_{4}=(\tilde{\beta}\tau^{k_{1}}w_{2})(\tilde{\beta}\tau^{k_{2}-1}w_{2}w_{4})\\
\tilde{\beta}\tau^{k_{1}-2}w_{2}w_{4}^{2}=d_{4}\tilde{\beta}\tau^{k_{1}}w_{2}\\
\begin{array}{cc}
\tilde{\beta}\tau^{-3}w_{2}w_{3}w_{4}^{2}=d_{3}d_{4} & \tilde{\beta}\tau^{k_{1}+k_{2}-2}w_{2}w_{3}w_{4}^{2}=d_{4}(\tilde{\beta}\tau^{k_{1}}w_{2})(\tilde{\beta}\tau^{k_{2}}w_{2})\end{array}
\end{array}
\end{equation}

\noindent As a more detailed example, we can consider why $\tilde{\beta}\tau^{-1}w_{2}w_{3}$
is equal to $d_{3}$. We know that $t(\tilde{\beta}\tau^{-1}w_{2}w_{3})$
and $t(d_{3})$ are both equal to $(\tilde{\beta}w_{2})^{-1}$ and
both have degree $(6,3)$. Our choice of notation gives $\mu(d_{3})=\tau^{-1}(w_{3}^{alg})^{2}$,
so Lemma \ref{lem: mod 2 reductions equal means originals equal}
and \cite{Por21} 7.3.1 give that $\mu(\tilde{\beta}\tau^{-1}w_{2}w_{3})=\tau^{-1}(w_{3}^{alg})^{2}$
and therefore $\tilde{\beta}\tau^{-1}w_{2}w_{3}=d_{3}$. In short,
the only integral generators that appear are $y_{2}$, $d_{2}$, $d_{3}$,
$d_{4}$, $\tilde{\beta}\tau^{k}w_{2}$ and $\tilde{\beta}\tau^{k-1}w_{2}w_{4}$.

It remains to determine how these six types of generators interact
with each other. First, recall the original relations from Field:

\begin{equation}
\begin{array}{ccc}
2d_{3}=0 & y_{2}d_{3}=0 & y_{2}^{2}-4d_{4}=0\end{array}
\end{equation}

\noindent Next, the commutative diagram construction shows that the
new classes are all $2$-torsion:

\noindent 
\begin{equation}
\begin{array}{cc}
2\tilde{\beta}\tau^{k}w_{2}=0 & 2\tilde{\beta}\tau^{k-1}w_{2}w_{4}=0\end{array}
\end{equation}

\noindent Lemma \ref{lem: mod 2 reductions equal means originals equal}
along with $\mu(\tilde{\beta}\tau^{k}w_{2})=\tau^{k}w_{3}^{alg}$
and $\mu(\tilde{\beta}\tau^{k}w_{2}w_{4})=\tau^{k}w_{3}^{alg}w_{4}^{alg}$
give the following:

\noindent 
\begin{equation}
\begin{array}{c}
(\tilde{\beta}\tau^{k_{1}}w_{2})(\tilde{\beta}\tau^{k_{2}}w_{2})=(\tilde{\beta}\tau^{k_{3}}w_{2})(\tilde{\beta}\tau^{k_{4}}w_{2})\iff k_{1}+k_{2}=k_{3}+k_{4}\\
(\tilde{\beta}\tau^{k_{1}-1}w_{2}w_{4})(\tilde{\beta}\tau^{k_{2}-1}w_{2}w_{4})=(\tilde{\beta}\tau^{k_{3}-1}w_{2}w_{4})(\tilde{\beta}\tau^{k_{4}-1}w_{2}w_{4})\iff k_{1}+k_{2}=k_{3}+k_{4}\\
(\tilde{\beta}\tau^{k_{1}}w_{2})(\tilde{\beta}\tau^{k_{2}-1}w_{2}w_{4})=(\tilde{\beta}\tau^{k_{3}}w_{2})(\tilde{\beta}\tau^{k_{4}-1}w_{2}w_{4})\iff k_{1}+k_{2}=k_{3}+k_{4}
\end{array}
\end{equation}
Finally, $d_{2}$ does not interact with any other generators and
interactions of new generators with $y_{2}$ are similar to those
between $y_{2}$ and $d_{3}$:

\noindent 
\begin{equation}
\begin{array}{cc}
y_{2}\tilde{\beta}\tau^{k}w_{2}=0 & y_{2}\tilde{\beta}\tau^{k-1}w_{2}w_{4}=0\end{array}
\end{equation}

These final equations are shown in the same way one shows that $y_{2}d_{3}=0$.
Yagita's presentation shows that $y_{0,2}\in H^{*,*'}(BSO_{4};\mathbb{Z}_{2})$
only appears as part of $\mathbb{Z}_{2}[\mu(d_{2}),\mu(d_{4})]\{y_{0,2}\}$.
This implies that $y_{0,2}$ multiplied by any generator other than
$\mu(d_{2})$ or $\mu(d_{4})$ is $0$ (\cite{HN18} Theorem 0.1).
Keeping in mind that $2\tilde{\beta}w_{2}=0$, $t(y_{2})=2\sqrt{p_{2}}$,
$t(d_{3})=(\tilde{\beta}w_{2})^{2}$, $t(\tilde{\beta}\tau^{k}w_{2})=\tilde{\beta}w_{2}$
and $t(\tilde{\beta}\tau^{k}w_{2})=\sqrt{p_{2}}\tilde{\beta}w_{2}$:

$$\xymatrix@!C=5pc{
y_{2}d_{3} \ar@{|->}[r]^-{\mu} \ar@{|->}[d]^{t} & y_{0,2}\mu(d_{3}) \ar@{|->}[d]^{t_{2}} & y_{2}\tilde{\beta}\tau^{k}w_{2} \ar@{|->}[r]^-{\mu} \ar@{|->}[d]^{t} & \tau^{k}y_{0,2}w_{3}^{alg} \ar@{|->}[d]^{t_{2}} & y_{2}\tilde{\beta}\tau^{k-1}w_{2}w_{4} \ar@{|->}[r]^-{\mu} \ar@{|->}[d]^{t} & \tau^{k-1}y_{0,2}w_{2}^{alg}w_{4}^{alg} \ar@{|->}[d]^{t_{2}} \\
0 \ar@{|->}[r]^-{\mu} & 0 & 0 \ar@{|->}[r]^-{\mu} & 0 & 0 \ar@{|->}[r]^-{\mu} & 0
}$$Thus, the upper right entry of each of these four diagrams must be
equal to $0$. We can demonstrate that the upper left corners must
all be $0$ as well using $y_{2}d_{3}$ as an example:
\begin{itemize}
\item $\mu_{M}(y_{2}d_{3}\otimes1)=y_{0,2}\mu(d_{3})$ and $\mu_{M}$ is
injective, so $y_{2}d_{3}\otimes1=0$
\item $y_{2}$ and $d_{3}$ are both not $2$-divisible, so therefore $y_{2}d_{3}=0$
\end{itemize}
\noindent The claim that $y_{2}\tilde{\beta}\tau^{k}w_{2}=0$ and
$y_{2}\tilde{\beta}\tau^{k-1}w_{2}w_{4}=0$ follows from the same
logic because $\tilde{\beta}\tau^{k}w_{2}$ and $\tilde{\beta}\tau^{k-1}w_{2}w_{4}$
are not $2$-divisible.
\noindent \begin{flushright}
$\boxempty$
\par\end{flushright}

\section{Conclusion}

The above represents a portion of the author's thesis work, the full
version of which can be found at \cite{Por21}. There are several
additional results presented there, most notably the analogous computation
for the integral motivic cohomology of $BG_{2}$ using the results
presented in this paper. A future goal is to also make the $BG_{2}$
part of this work a standalone paper.

\end{document}